\documentclass[12pt,a4paper]{article}
\usepackage{amsmath,euscript}
\usepackage{amssymb}
\usepackage{amsthm}
\usepackage{enumerate}
\usepackage{amsfonts}
\usepackage{comment}
\usepackage{amscd}

\newtheorem{theorem}{Theorem}[section]

\newtheorem{proposition}{Proposition}[section]
\newtheorem{lemma}[theorem]{ Lemma}%[section]
\newtheorem{Proposition}{\bf Proposition}[section]
\newtheorem{question}{\bf Question}[section]
\newtheorem{corollary}[theorem]{Corollary}%[section]

\title{A note on Nil and Jacobson radicals in graded rings}

\author{Agata Smoktunowicz}

\date{ }

\begin{document}

\maketitle

\begin{abstract}  It was shown by Bergman that the Jacobson radical of a $Z$-graded ring is homogeneous. This paper shows that the analogous result holds for nil rings, namely, that the nil radical of a $Z$-graded ring is homogeneous.

 It is obvious that a subring of a nil ring is nil, but generally a subring of a Jacobson radical ring need not be a Jacobson radical ring. In this paper it is shown  that  every subring which is generated by homogeneous elements in a graded Jacobson radical ring is always a Jacobson radical ring. It is also observed that a ring whose all subrings are Jacobson radical rings is nil. Some new results on graded-nil rings are also obtained.
 % We finish our paper by answering in the negative the question of Shestakov, by showing that there %are locally nilpotent rings such that the rings of skew polynomials over these rings are not nil (not %Jacobson radical).
\end{abstract}
  It was shown by Bergman that the Jacobson radical of a $Z$-graded ring is homogeneous \cite{b}.
 During the conference  'Radicals in Rings and Related Topics' held in Warsaw in 2009, the question was posed whether the same is true for the nil radical. Recall that the nil radical of a ring $R$ is the largest nil ideal in $R$. It is known that nil and Jacobson radical in polynomial rings are homogeneous. Indeed, it was proved by Amitsur that for every ring $R$ the Jacobson radical
$J(R[x])$ of the polynomial ring $R[x]$ in an indeterminate $x$ over $R$ is equal to $I[x]$
for a nil ideal $I$ of $R$ \cite{a2}.
The original proof of this theorem by Amitsur also shows that for every ring $R$ the nil radical $N(R[x])$ of the polynomial
ring $R[x]$ in an indeterminate $x$ over $R$ is equal to $I[x]$ for some nil ideal $I$ of $R$. In a more general setting, in a series of papers  Jespers, Keralev,  Krempa, Okni\'nski and Puczy\l{}owski proved that for many groups the Jacobson radical of a $G$-graded ring is homogeneous \cite{cj, cr, j1, j2, ep, jo, ko}. In this paper we show that nil radicals in graded rings have an analogous property.
  The main result of this paper is as follows.
 \begin{theorem}
  The nil radical of a $Z$-graded ring is homogeneous.
 \end{theorem}
  Recall that a $Z$-graded ring is a ring graded by the additive group of integers, and that an ideal $I$ in a graded ring is homogeneous if it is generated by homogeneous elements. Recall also that a free abelian group is always a direct product of the integers $Z$. Hence, from Theorem $0.1$ we get the following corollary.
   \begin{corollary}
   Let $G$ be a free abelian group, then the nil radical of a $G$-graded ring is homogeneous.
   \end{corollary}
 %   Recall also that an element in a graded ring $R=\oplus_{i=0,1, \ldots }R_{i}$ is homogeneous if it belongs to $R_{i}$ for some $i$. Naturally an ideal $I$ in $R$ is % homogeneous if and only if it is generated by homogeneous elements of $R$.

  It is important to ascertain  when nil and Jacobson radicals coincide. It is known that nil rings are Jacobson radical, and on the other hand Amitsur's theorem assures that the Jacobson radical of a finitely generated algebra over an uncountable field is nil \cite{a1}.
  This is important partly because nil and Jacobson radical have different properties, for example matrix rings over Jacobson radical rings are Jacobson radical, but  whether rings of matrices over nil rings are nil is a tantalizing open question. In fact it was shown by Krempa \cite{k} in 1972 that this question is equivalent to  Koethe's famous conjecture of 1930: {\em If R has a non-zero nil one sided ideal, does it follow that R has a non-zero nil two-sided ideal}? This question remains unanswered. On the other hand matrix rings over graded-nil rings need not be graded-nil  (graded-nil rings are graded rings in which all homogeneous elements are nilpotent) \cite{s}. Nil rings also possess some properties which Jacobson radical rings do not possess, for example a subring of a nil ring is nil, but a subring of a Jacobson radical ring need not be Jacobson radical. Recall that a ring $R$ is called Jacobson radical if the Jacobson radical of $R$ is $R$.
   In this paper we observe that a homogeneous subring of the Jacobson radical in a graded (by the semigroup of positive integers) ring is Jacobson radical; that is, we prove the following result:
    \begin{Proposition}
    Let $R$ be a Jacobson radical ring graded by the semigroup of positive integers. Then every subring of $R$ which is generated by homogeneous elements is a Jacobson radical ring.
    \end{Proposition}
     Notice that, if $S$ is a subring generated by a single element $a$, then  $S$ is a Jacobson radical ring if and only if $a$ is nilpotent.
  Therefore a Jacobson radical ring whose all subrings are Jacobson radical rings is nil.

 It would also be interesting to know if more general results hold; therefore we  pose the following two  questions:
 \begin{question} For which groups $G$ is the nil radical of a $G$-graded ring homogeneous?
 \end{question}
\begin{question} For which semigroups $G$ are homogeneous subrings of $G$-graded Jacobson radical rings also Jacobson radical rings?
\end{question}
 For interesting results on $G$-graded rings see \cite{cr}. Note that radicals of group rings and semigroup rings were studied even before Bergman's preprint \cite{k2} (see also \cite{ks}). Results from \cite{jkp} are related to the open questions $0.1$, $0.2$.

    \section{ Nil rings}

 The purpose of this section is to prove Theorem $0.1$. We begin by proving the following lemma.

 \begin{lemma}
 Let $R=\oplus _{i\in Z}R_{i}$ be a graded ring and $N$ be the nil radical of $R$. If $N$ contains no nonzero homogeneous elements then $N=0$.
 \end{lemma}

{\bf Proof.}
  We will prove that for every $k,n\in Z$, $n\geq 0$, if $a=a_{k}+\ldots +a_{k+n}\in N$ with $a_{i}\in R_{i}$
    then $a_{k}=a_{k+1}=\ldots =a_{k+n}=0$. We proceed by induction on $n$.

   If $n=0$ we get that $a=a_{k}\in N$, and since $N$ does not  contain any nonzero homogeneous elements and $a_{k}$ is homogeneous we get $a=a_{k}=0$.

  Let now $n>0$, and suppose that the result is true for all numbers smaller than $n$. We need to show that it holds for $n$.
   Recall that the prime radical (the intersection of all prime ideals) in a $Z$-graded ring is homogeneous and contained in $N$ \cite{l}. Therefore the prime radical of $R$ is zero (because by assumption $N$ contains no nonzero homogeneous elements). It
 follows that $R$ is a semiprime ring.

 For each homogeneous $r\in R$, $a_{k+n}ra-ara_{k+n}\in N$, since $a=a_{k}+\ldots +a_{k+n}\in N$.
  It follows that $a_{k}'+a_{k+1}'+\ldots +a_{k+n-1}'\in N$ where $a_{i}'=a_{n+k}ra_{i}-a_{i}ra_{n+k}$.
   By the inductive assumption $a_{k}'=a_{k+1}'=\ldots =a_{k+n-1}'$.
    It follows that $a_{k+n}ra_{i}-a_{i}ra_{k+n}=0$ for all homogeneous $r\in R$, so that $a_{k+n}ra-ara_{k+n}=0$ for all $r\in R$.

 To get the desired result it suffices to prove that $a_{k+n}=0$ because then  by the inductive assumption $a_{k}=\ldots =a_{k+n-1}=0$, so $a=a_{k}+\ldots +a_{k+n-1}=0$.

   Suppose on the contrary that $a_{k+n}\neq 0$. As $a_{k+n}$ is homogeneous, it follows that $a_{k+n}\notin N$.  That means that there are  $p_{i},q_{i}\in R$, such that the element $\sum_{i} p_{i}a_{n+k}q_{i}$ is not nilpotent.
    Notice that $a\in N$ implies $(\sum_{i} p_{i}aq_{i})^{d}=0$ for some $d$. If $(\sum_{i} p_{i}a_{n+k}q_{i})^{d}=0$ then we have a contradiction and the proof is finished. Hence, we assume that $(\sum_{i} p_{i}a_{n+k}q_{i})^{d}\neq 0$.
   If $(\sum_{i} p_{i}a_{n+k}q_{i})^{d}Ra\neq 0$ then since $R$ is semiprime
   $$((\sum_{i} p_{i}a_{n+k}q_{i})^{d}RaR)^{d}\neq 0$$
    (this follows because in a semiprime ring, if $a\neq 0$, then $(aR)^{m}\neq 0$ for all $m$). It follows from $a_{n+k}ra=ara_{n+k}$ that $(\sum_{i} p_{i}a_{n+k}q_{i})RaR=(\sum_{i} p_{i}aq_{i})Ra_{n+k}R\subseteq (\sum_{i} p_{i}aq_{i})R$. By applying this argument again we get that
   $(\sum_{i} p_{i}a_{n+k}q_{i})^{2}RaRaR\subseteq (\sum_{i} p_{i}aq_{i})^{2}R$. Continuing in this way we get
     $(\sum_{i} p_{i}a_{n+k}q_{i})^{d}(RaR)^{d} \subseteq(\sum_{i} p_{i}aq_{i})^{d}R=0$. Observe that $((\sum_{i} p_{i}a_{n+k}q_{i})^{d}RaR)^{d}\subseteq  (\sum_{i} p_{i}a_{n+k}q_{i})^{d}(RaR)^{d}=0$ a contradiction.
     Therefore $$(\sum_{i} p_{i}a_{n+k}q_{i})^{d}Ra=0.$$

 We will now show that $(\sum_{i} p_{i}a_{n+k}q_{i})^{d}Ra_{n+k}=0$.
  Recall that $$(\sum_{i} p_{i}a_{n+k}q_{i})^{d}Ra=0.$$ Write $(\sum_{i} p_{i}a_{n+k}q_{i})^{d}=e_{t}+\ldots +e_{j}$ with each  $e_{i}\in R_{i}$ for some $t,j\in Z$, $t\leq j$ and $e_{j}\neq 0$.
  Recall that $a=a_{k}+\ldots +a_{k+n}$. Let $r$ be a homogeneous element of $R$. Now $(e_{t}+\ldots +e_{j})r(a_{k}+\ldots +a_{k+n})=0$ implies $e_{j}ra_{k+n}=0$ (by comparing elements of given degree as the ring is graded). Moreover, for each $v<j$, $e_{v}ra_{n+k}=-\sum_{i>v} e_{i}ra_{n+k+v-i}$.
    By applying this last equation repeatedly to elements from  $(e_{t}+\ldots +e_{j})(Ra_{n+k})^{j}$ we get $(e_{t}+\ldots +e_{j})(Ra_{n+k})^{j}\subseteq e_{j}R$.
     Therefore, $(e_{t}+\ldots +e_{j})(Ra_{n+k})^{j+1}\subseteq e_{j}Ra_{n+k}$. Recall that $e_{j}Ra_{n+k}= 0$,
     and so $(e_{t}+\ldots +e_{j})(Ra_{n+k})^{j+1}=0$. Observe now that
     $((e_{t}+\ldots +e_{j})Ra_{n+k}R)^{j+1}\subseteq (e_{t}+\ldots +e_{j})(Ra_{n+k}R)^{j+1}\subseteq (e_{t}+\ldots +e_{j})(Ra_{n+k})^{j+1}R=0$.
     Recall that $R$ is a semiprime ring, hence  $((e_{t}+\ldots +e_{j})Ra_{n+k}R)^{j+1}=0$ implies $(e_{t}+\ldots +e_{j})Ra_{n+k}=0$. Because
     $(\sum_{i} p_{i}a_{n+k}q_{i})^{d}=e_{t}+\ldots +e_{j}$  we have $(\sum_{i} p_{i}a_{n+k}q_{i})^{d}Ra_{n+k}=0$, and so $(\sum_{i} p_{i}a_{n+k}q_{i})^{d}Ra_{n+k}R=0$.

  Notice that $\sum_{i} p_{i}a_{n+k}q_{i}\subseteq Ra_{n+k}R$,  implies
 $(\sum_{i} p_{i}a_{n+k}q_{i})^{d}R(\sum_{i} p_{i}a_{n+k}q_{i})^{d}=0$ and since $R$ is semiprime we get
 $(\sum_{i} p_{i}a_{n+k}q_{i})^{d}=0$, a contradiction.
 This finishes the proof.

      We will now prove Theorem $0.1$.

     {\bf Proof of Theorem $0.1$}. Let $M$ be the sum of all
      nil ideals in $R$ which are generated by graded elements. Let $<a>$ denote the ideal generated by the element $a$. Hence
     $M=\sum_{a\in W} <a>$ where $W$ is the set of all homogeneous elements of $R$ which are in the nil radical of $R$.
       Clearly $M$ is a nil ideal, and $M$ is homogeneous. Consider ring $R/M$. Let $a\in R/M$ be a homogeneous element which is in the nil ideal of $R/M$. We claim that $a=0$.
        Since $a$ is homogeneous, $a=r+M$ for some homogeneous element $r\in R$. Since $a$ is in the nil radical of $R/M$ and $M$ is nil, it follows that $r$ is in the nil radical of $R$, and hence $r\in M$, so $a=r+M=M$ is zero in $R/M$. It follows that the nil radical of $R/M$ does not contain any non-zero homogeneous elements, so by the previous lemma it is zero. It follows that the nil radical of $R$ equals $M$ (if it was bigger than $M$ then $R/M$ would contain a nonzero nil ideal).

We will now prove Corollary $0.1$.

{\bf Proof of Corollary $0.1.$} Free abelian groups are direct products of copies of $Z$. Therefore we can write $G=\oplus_{i\in F}Z_{i}$ with each $Z_{i}$ isomorphic to $Z$, and $F$ being a set. Fix an index $i$. Our algebra is graded by the group $G$, and so it is also graded by the group $Z_{i}$. By Theorem $0.1$ the nil radical of our algebra is $Z_{i}$-- graded. Notice that it holds for every index $i$, hence the nil radical is $G$-graded. The proof is finished.

      \section{Graded Jacobson radical rings}

  All rings in this section are graded by the additive semigroup of positive integers. The purpose of this section is to prove Proposition $0.1$.

  {\bf Proof of Proposition $0.1.$} Let $R=\oplus_{i=1}^{\infty }R_{i}$ be a Jacobson radical ring graded by the additive semigroup of positive integers. Denote
   $R^{*}=\{\sum_{i=1}^{\infty}r_{i}:r_{i}\in R_{i}\}$ (by analogy to power series rings). Observe that the ring $R$ can be viewed as a subring of $R^{*}$. Notice also that since $R$ is graded then $R^{*}$ has well defined multiplication and addition, and satisfies all axioms of associative ring. It is known \cite{l} that a ring $R$ is Jacobson radical, if and only if for every $a\in R$ there is $a'\in R$ such that $a+a'+aa'=a+a'+a'a=0$, such an element $a'$ is called a quasi-inverse of $a$ \cite{l}. It is known that in a Jacobson radical ring the quasi-inverse of each element is uniquely determined \cite{l}. Observe that ring $R^{*}$ is Jacobson radical and the quasi-inverse of an element $a\in R^{*}$ is $a'=\sum_{i=1}^{\infty }(-1)^{i}a^{i}$ which is a well defined element of $R^{*}$. Indeed $a+a'+aa'=0$.
  Let $a\in R$, then $a+a''+aa''=0$ for some $a''\in R$. The quasi-inverse of an element $a$ in  Jacobson radical ring $R^{*}$ is uniquely determined, hence  $a''=a'=\sum_{i=1}^{\infty }(-1)^{i}a^{i}$, and so the series $\sum_{i=1}^{\infty }(-1)^{i}a^{i}$ has almost all homogeneous components equal to zero.

  Let $S$ be a subring of $R$ generated by homogeneous elements.
 Let $a\in S$, and let $a=a_{1}+a_{2}+\ldots +a_{n}$ with $a_{i}\in R_{i}$. Notice that $a_{i}\in S$ for all $i\leq n$, since $S$ is homogeneous.  Recall that $R$ is Jacobson radical, hence $a$ has a quasi-inverse $a'$ in $R$,
 $a'=\sum_{i=1}^{k }r_{i}$ for some $k$ and some $r_{i}\in R_{i}$. Let $a''=\sum_{i=1}^{\infty }(-1)^{i}a^{i}$ be the quasi-inverse
 of $a$ in $R^{*}$. We know that $a'=a''$, so that $\sum_{i=1}^{\infty }(-1)^{i}a^{i}=\sum_{i=1}^{k }r_{i}$.
    Notice that all $r_{i}\in S$ because  all $a_{i}\in S$. Recall that almost all $r_{i}$ are zero. It follows that the quasi-inverse of $a$ in $R$ is $a'=\sum_{i=0}^{d}r_{i}$ for some $d$, with all $r_{i}\in S$, and so $a'\in S$. Hence $S$ is a Jacobson radical ring.

{\bf Remark.} In the 1970's Krempa introduced the idea of embedding polynomial ring $R[x]$ into a power series ring $R\{x\}$ to study the Jacobson radical of polynomial rings. This idea inspired the author to invent Proposition $0.1$.

 Observe that polynomial ring $R[x,x^{-1}]$ is $Z$-graded, with $x$ having gradation one and elements of $R$ having gradation $0$. It follows that nil and Jacobson radical of  ring $R[x,x^{-1}]$ are homogeneous. We pose the following open question:
 \begin{question}
 Let $R$ be a nil ring. Is $R[x,x^{-1}]$ a Jacobson radical ring?
 \end{question}
 Note that it was proved by Krempa that Koethe's conjecture is equivalent to the following assertion: {For every nil ring $R$, the polynomial ring $R[x]$ is Jacobson radical} \cite{k}. We don't know if Question $2.1$ is equivalent to Koethe's conjecture. We also don't know if there is a ring $R$ such that $R[x]$ is Jacobson radical but $R[x,x^{-1}]$ is not Jacobson radical.

      \section{Graded-nil rings and the Brown-McCoy radical}

  Recall that the prime radical in a ring $R$ is the intersection of all prime ideals in $R$. Similarly, in a ring with an identity element, the Brown-McCoy radical is the intersection of all maximal ideals. In general, given ring $R$ the Brown-McCoy radical $U(R)$ of $R$ is the intersection of all ideals $I$ of $R$ such that $R/I$ is a simple ring with an identity element. In particular a ring is Brown-McCoy radical if and only if it cannot be homomorphically mapped onto a simple ring possessing an  identity element. It is obvious that every Jacobson radical ring is Brown-McCoy radical. We will show that for a ring graded by the additive semigroup of positive integers the following slightly stronger result holds.
  \begin{proposition}
  Let $R$ be a ring graded by the additive semigroup of positive integers. If $R$ is graded-nil then $R$ is Brown-McCoy radical.
  \end{proposition}
   Recall that a graded ring is graded-nil if all homogeneous elements in $R$ are nilpotent. It is known that graded-nil rings even over uncountable fields need not be Jacobson radical \cite{s}. On the other hand every graded Jacobson radical ring is graded-nil.

  {\bf Proof of Proposition $3.1$} Suppose that $R$ is not Brown-McCoy radical. Then $R/I$ is a simple ring with an identity element, for some proper ideal $I$ in $R$.
  Recall that an ideal $P$ of a ring (without unity) is (right) primitive if and only if there exists a modular maximal right ideal $Q$ of $R$ such that $P$ is the maximal two-sided ideal contained in $Q$. Recall also that
 a right ideal $Q$ in a ring $R$ is modular if and only if there exists element $b\in R$ such that $a-ba\in Q$ for every $a\in R$.

    Observe that $I$ is a maximal ideal and since $R/I$ has an identity element $I$ is also a primitive ideal. To see that
    $I$ is primitive proceed as follows: let $e\in R$ be such that $e+I$ is the identity element in $R/I$. Then $ea-a\in I$  for all $a\in R$ hence $I$ is a modular right ideal in $R$. Notice that $e\notin I$ as otherwise $I=R$. Let $Q$ be a maximal right ideal in $R$ which contains $I$ and does not contain $e$ (it exists by Zorn's lemma).
    Then $Q$ is a maximal modular right ideal in $R$, and since $R/I$ is simple then $I$ the largest two-sided ideal contained in $Q$. It follows that $I$ is a primitive ideal in $R$.

  It was shown in \cite{s2} that primitive ideals in graded nil rings are homogeneous. It follows that $R/I$ is graded. However, a graded (by positive integers) ring cannot have an identity element, as for large $n$, $e^{n}$ has bigger degree than $e$, so there are no nonzero elements satisfying $e^{n}=e$.

  {\bf Remark} It is not known if Proposition $3.1$ also holds for $Z$ graded algebras, or whether primitive ideals in $Z$ graded rings which are graded-nil are homogeneous. We also don't know whether polynomial ring over a graded Jacobson radical ring $R=\oplus_{i=1}^{\infty }R_{i}$ is always Brown-McCoy radical.

 {\bf Acknowledgements.} The author is very grateful to Tom Lenagan, Jan Okni\'nski and Jan Krempa for many helpful suggestions which improved the paper.

\vspace{3ex}

\begin{minipage}{1.00\linewidth}

\noindent
Agata Smoktunowicz: \\
School of Mathematics, \\
The University of Edinburgh,\\
James Clerk Maxwell Building, \\
The King's Buildings, Mayfield Road\\
EH9 3JZ, Edinburgh. \\

\end{minipage}

\end{document}